\documentclass[a4,11pt,reqno]{amsart}%article}
\usepackage[english]{babel}
\usepackage[latin1]{inputenc}
\usepackage{graphics}
\usepackage{graphicx}
\usepackage{ulem}
\usepackage{hhline}
\usepackage{dsfont}
\usepackage{mathrsfs}
\usepackage{fancyhdr}
\usepackage{amsmath,amssymb}
\usepackage{rotating}
\usepackage[latin1]{inputenc}
\usepackage[T1]{fontenc}
\usepackage{fancybox}
\usepackage{color}
\usepackage{colortbl}
\usepackage{setspace}
\usepackage{enumerate}
\usepackage[nice]{nicefrac}
\usepackage{amsthm}
\usepackage{wasysym}
\usepackage{epsfig}
\usepackage{multicol}
\usepackage{pifont}
\usepackage[french]{minitoc}
\usepackage{latexsym,amsfonts}
\usepackage{varioref}
\usepackage{textcomp}
\usepackage{lmodern}
\usepackage{mathpazo}
\usepackage{euscript}
\usepackage[pdftex]{hyperref}

\numberwithin{equation}{section} \makeatletter\@addtoreset{equation}{section}

\DeclareMathSymbol{\subsetneqq}{\mathbin}{AMSb}{36}

\newcommand{\set}[1]{{\left\{{#1}\right\}}}

\newcommand{\scal}[1]{{\left\langle{#1}\right\rangle}}
\newcommand{\C}{\mathbb C}
\newcommand{\R}{\mathbb R}

\newtheorem {theorem}{Theorem}[section]

\newtheorem {lemma}[theorem]{Lemma}
\newtheorem {proposition}[theorem]{Proposition}
\newtheorem {remark}[theorem]{Remark}

\newtheorem {corollary}[theorem]{Corollary}

\begin{document}

\title[A characterization of planar mixed automorphic forms]{A characterization of planar mixed automorphic forms}
\author{A. Ghanmi}
%\address{Department of Mathematics,\newline  Faculty of  Sciences, P.O. Box 1014\newline
%   Mohammed V University,  Agdal, \newline 10 000 Rabat - Morocco  }
%   \email{allalghanmi@gmail.com}
%\keywords{11F12, 11F55, 32N10}

\date{ %Novembre 27, 2009
\today
\newline {\it E-mail address}: allalghanmi@gmail.com}
\maketitle

\vspace*{-.8cm}
\begin{center}
{\it Department of Mathematics,  Faculty of  Sciences, P.O. Box 1014, \\
   Mohammed V University,  Agdal,  10000 Rabat - Morocco}
   \end{center}
\begin{abstract} \quad
We characterize the space of the so-called planar mixed automorphic forms of
type $(\nu,\mu)$ with respect to an equivariant pair $(\rho,\tau)$ as the image, by an appropriate transform, of the usual (Landau)
 automorphic forms involving special multiplier.
\end{abstract}

 \section{Introduction}

 The notion of mixed automorphic forms was introduced by Stiller \cite{Stiller84} and
 extensively studied by M.H. Lee (see \cite{LeeLectNotes04} and the references therein).
 They appears essentially in the context of number theory and algebraic geometry
 and arise naturally as holomorphic forms on elliptic varieties \cite{HuntMeyer85}. Mixed automorphic forms
   include classical ones as a special case and % (take for example $\mu=0$ or let $\rho$ and $\tau$ be the identity maps).
 non trivial examples of them has being constructed in \cite{Choie97,Lee00}.
 In this paper, we deal with the space planar mixed automorphic forms and we show that it can be connected  to the space of Landau automorphic forms \cite{GI-JMP08}, by an explicit and special transform (Theorem \ref{Thm2}).

 Let $\C$ be the complex plane endowed with its usual hermitian scalar product $\scal{z,w}=z\bar w$, and $\mathbf{T}$ be
 the unitary group, $\mathbf{T}=\set{\lambda\in \C; |\lambda|=1}$. Consider the semidirect product group $G=\mathbf{T}\rtimes \C$
  operating on $\C$ by the holomorphic mappings $g\cdot z = a z +b$ for $g=(a,b)\in G$. By equivariant pair $(\rho,\tau)$,
%\cite{Kuga64,SatakeBook80,HuntMeyer85,Abdulali95,LeeLectNotes04}.
we mean that  $\rho$ is a $G$-endomorphism  and $\tau : \C \to \C$ a compatible mapping  such that
\begin{equation}
\label{EquivCond} \tau(g.z)= \rho(g)\cdot\tau(z); \qquad g\in G, \quad z\in\C.
\end{equation}
 Associated to such $(\rho,\tau)$ and given uniform lattice $\Gamma$ in $\C$,
we consider the vector space $\mathcal{M}^{\nu,\mu}_{\tau}(\C)$
of $\Gamma$-mixed automorphic forms of type
$(\nu,\mu)$.  They are smooth complex-valued functions $F$ on $\C$
satisfying the functional equation
 \begin{equation}\label{FunctionalEq1}
 F(\gamma\cdot z) =  j^\nu(\gamma, z)j^\mu(\rho(\gamma), \tau(z)) F(z),
\end{equation}
where $\nu, \mu$ are non negative real numbers and $j^\alpha$; $\alpha\in \R$, is defined %to be the automorphic factor given
by
\begin{equation}\label{AutomFactorExplicit}
 j^\alpha(g,z)=e^{2i\alpha  \Im \scal{z,g^{-1}\cdot 0}}.
\end{equation}
Here and elsewhere $\Im z$ denotes the imaginary part of the complex number $z$. We assert that the space
 $\mathcal{M}^{\nu,\mu}_{\tau}(\C)$ is isomorphic to the space of automorphic forms $F\in \mathcal{C}^\infty(\C)$, i.e.,
 such that
 $$ F(z+\gamma)= \chi_\tau(\gamma)j^{B^{\nu,\mu}_{\tau}}(\gamma, z)  F(z) ;\qquad z\in \C, \gamma\in \Gamma
  %:=\mathcal{F}^{B^{\nu,\mu}_{\tau}}_{\Gamma,\chi_{\tau}}
.$$
The pseudo-character $\chi_\tau$ is defined on $\Gamma$ by
 $$\chi_\tau(\gamma)=\exp\bigg({2i\varphi^{\nu,\mu}_\tau(\gamma)-2i\mu  \Im \scal{\tau(0),\rho(\gamma)^{-1}\cdot 0}}\bigg),$$
where the involved function $\varphi^{\nu,\mu}_\tau$ satisfies a first order differential equation as in Proposition \ref{1stdiffEq} below.

Our main result (Theorem \ref{Thm2}) is stated and proved in Section 3. To do this, we have  to ensure first the nontrivially of the space $\mathcal{M}^{\nu,\mu}_\tau(\C)$ and to introduce properly the function $\varphi^{\nu,\mu}_\tau$ (Section 2). The crucial point in the proof of Theorem \ref{Thm2} is to observe that the weight $B^{\nu,\mu}_\tau$ given by  $B^{\nu,\mu}_\tau = \nu +\mu(|\frac{\partial \tau}{\partial z}|^2- |\frac{\partial \tau}{\partial\bar z}|^2)$ is indeed a real constant independent of the complex variable $z$.  As immediate application of the obtained characterization, one can deduce easily some concrete spectral properties of an appropriate invariant Laplacian
  acting on $\mathcal{M}^{\nu,\mu}_{\tau}(\C)$ (see \cite{GhanmiArxiv09} for more details).
% Mainly, we see that the dimensional formula of the eigenspace $$\mathcal{E}^{\nu,\mu}_{\tau;k}:= \set{F\in\mathcal{M}^{\nu,\mu}_{\tau}(\C); \quad \La F = B^{\nu,\mu}_\tau (2k+1)F} $$ is given explicitly by
%$\dim \mathcal{E}^{\nu,\mu}_{\tau;k} = (2B^{\nu,\mu}_\tau/\pi) Area(\C/\Gamma).$

 \section{On the space $\mathcal{M}^{\nu,\mu}_{\tau}(\C)$ and the function $\varphi^{\nu,\mu}_\tau$}

For given real numbers $\nu,\mu>0$ and given equivariant pair $(\rho,\tau)$, we define
 $J^{\nu,\mu}_{\rho,\tau}$ as a complex valued mapping  $ G\times \C$ by
\begin{equation}\label{AutomFactor}
J^{\nu,\mu}_{\rho,\tau} (g,z):= j^\nu(g, z) j^\mu(\rho(g), \tau(z))
\end{equation}
and we perform the vector space of mixed automorphic forms of type $(\nu,\mu)$,
\begin{equation}\label{MixedForms}
\mathcal{M}^{\nu,\mu}_{\tau}(\C):= \set{F: \C \stackrel{\mathcal{C}^\infty}{\longrightarrow} \C;  \quad  F(\gamma \cdot z) = J^{\nu,\mu}_{\rho,\tau}
(g,z) F(z) , \quad \gamma \in \Gamma, z\in \C}.
\end{equation}

Then, one can check the following

\begin{proposition}\label{Prop1}
i) Let $\phi^{\nu,\mu}_\rho(g,g')$ be the real-valued function defined on $G\times G$ by
\begin{equation}\label{PhaseFactor}
\phi^{\nu,\mu}_\rho(g,g'):= \Im  \Big(\nu \scal{g^{-1}\cdot
0,g'\cdot 0} + \mu \scal{\rho({g^{-1}})\cdot
0,\rho(g')\cdot 0} \Big).
\end{equation} Then, the mapping $J^{\nu,\mu}_{\rho,\tau}$ satisfies the chain rule
\begin{align}
J^{\nu,\mu}_{\rho,\tau} (gg',z) = e^{2i\phi^{\nu,\mu}_\rho(g,g')}
% e^{i \nu \Im \scal{g^{-1}\cdot 0,g'\cdot 0}} e^{i \mu \Im \scal{\rho({g^{-1}})\cdot 0,\rho(g')\cdot 0}}
 J^{\nu,\mu}_{\rho,\tau} (g,g'\cdot z)
J^{\nu,\mu}_{\rho,\tau} (g',z).  \label{ChainRule}
\end{align}

ii) The functional space $\mathcal{M}^{\nu,\mu}_{\tau}(\C)$ is nontrivial if and only if the real-valued function $(1/\pi)\phi^{\nu,\mu}_\rho$ in \eqref{PhaseFactor} takes integer values on $\Gamma\times \Gamma$.

\end{proposition}

\begin{remark} According to i) of Proposition \ref{Prop1}, the unitary transformations %$\mathcal{T}^{\nu,\mu}_g$; $g\in G$,
\begin{equation}\label{Transformation}
 [\mathcal{T}^{\nu,\mu}_gf](z):= \overline{J^{\nu,\mu}_{\rho,\tau}
(g,z)} f(g\cdot z) ,
\end{equation}
 for varying $g\in G$, define then a projective representation of the group $G$ on the space of $\mathcal{C}^\infty$
functions on $\C$. While the assertion ii) shows that $\mathcal{M}^{\nu,\mu}_{\tau}(\C)$
can be realized as the space of cross sections on a line bundle over the complex torus $\C/\Gamma$.
\end{remark}

\begin{proof} For every $g,g'\in G$ and $z\in \C$, we have
\begin{align*}
J^{\nu,\mu}_{\rho,\tau} (gg',z) = j^{\nu} (gg',z)j^{\mu}(\rho(gg'),\tau(z)) = j^{\nu} (gg',z)j^{\mu}(\rho(g)\rho(g'),\tau(z)).
\end{align*}
%for $\rho$ being a groups homomorphism.
Next, one can sees that the automorphic factor $j^{\alpha}(\cdot,\cdot)$ satisfies
$$ j^{\alpha}(hh',w) =  e^{2i\alpha \Im  \scal{h^{-1}\cdot
0,h'\cdot 0}} j^{\alpha}(h, h'\cdot w) j^{\alpha}(h',w); \qquad h,h'\in G, w\in\C.$$ This gives rise to
\begin{align*}
J^{\nu,\mu}_{\rho,\tau} (gg',z)
%&= e^{2i\nu \Im  \scal{g^{-1}\cdot
%0,g'\cdot 0}} j^{\nu}(g, g'\cdot z) j^{\nu}(g',z) e^{2i\mu \Im  \scal{\rho(g)^{-1}\cdot
%0,\rho(g')\cdot 0}} j^{\mu}(\rho(g), \rho(g')\cdot \tau(z)) j^{\mu}(\rho(g'),\tau(z))\\&
= e^{2i\phi^{\nu,\mu}_\rho(g,g')} j^{\nu}(g, g'\cdot z)  j^{\mu}(\rho(g), \rho(g')\cdot \tau(z)) J^{\nu,\mu}_{\rho,\tau} (g',z).
\end{align*}
Finally, \eqref{ChainRule} follows by making use of the equivariant condition $\rho(g')\cdot \tau(z)=\tau(g'\cdot z)$.

 The proof of ii) can be handled in a similar way as in \cite{GI-JMP08} making use of \eqref{ChainRule} combined with the equivariant condition (\ref{EquivCond}).
\end{proof}

%\begin{proof}
In order to prove the main result of this paper, we need to introduce the function $\varphi^{\nu,\mu}_\tau$.
 \begin{proposition}\label{1stdiffEq} The first order differential equation
\begin{eqnarray}\label{Existence1stDiffEqProp} \frac{\partial \widetilde{\varphi^{\nu,\mu}_\tau}}{\partial \bar z} = -i\mu\Big(( \tau \frac{\partial \bar\tau}{\partial \bar z}- \bar \tau\frac{\partial \tau}{\partial\bar z}) - (|\frac{\partial \tau}{\partial z}|^2- |\frac{\partial \tau}{\partial\bar z}|^2) z\Big)\end{eqnarray}
admits a solution $\widetilde{\varphi^{\nu,\mu}_\tau}: \C \to \C$ such that $\Im\widetilde{\varphi^{\nu,\mu}_\tau}$ is constant.
 \end{proposition}

\begin{proof} By writing the $G$-endomorphism $\rho: G\to G=\mathbf{T}\rtimes \C$ as $\rho(g)=[\chi(g),\psi(g)]$,
% we see that the involved maps $\chi: G\to \mathbf{T}$ and $\psi: G \to \C$ satisfy
%$$ \chi(gg') =\chi(g)\chi(g') \quad \mbox{and} \quad \psi(gg')=\psi(g) +\chi(g)\psi(g').$$ Now, by
and differentiating the equivariant condition $ \tau(g\cdot z) =\rho(g)\cdot \tau(z)= \chi(g)\tau(z) + \psi(g),$ it follows
$$ (\partial\frac{g\cdot z}{\partial z})\frac{\partial \tau}{\partial z}({g\cdot z})= \chi(g)\frac{\partial \tau}{\partial z}(z) \quad \mbox{and} \quad
   (\frac{\partial\overline{g\cdot z}}{\partial \bar z})\frac{\partial \tau}{\partial \bar z}({g\cdot z})= \chi(g)\frac{\partial \tau}{\partial \bar z}(z).$$
   Hence, for $({\partial g\cdot z})/({\partial z})$ and $\chi(g) $ being in $\mathbf{T}$, we deduce that
   \begin{align*}  B^{\nu,\mu}_\tau(g\cdot z) & = \nu +\mu\bigg(\bigg|\frac{\partial \tau}{\partial z}({g\cdot z})\bigg|^2- \bigg|\frac{\partial \tau}{\partial\bar z}({g\cdot z})\bigg|^2\bigg) \\ & = \nu +\mu\bigg(\bigg|\frac{\partial \tau}{\partial z}({z})\bigg|^2- \bigg|\frac{\partial \tau}{\partial\bar z}({z})\bigg|^2\bigg) = B^{\nu,\mu}_\tau(z),
 \label{Super0}\end{align*}
 and therefore $z \mapsto B^{\nu,\mu}_\tau(z)$ is a real valued constant function (since the only $G$-invariant functions on $\C$ are the constants).
 Now, by considering the differential differential $1$-form $$\theta^{\nu,\mu}_\tau(z) :=i \Big\{\nu(\bar z dz -  zd\bar z) +\mu(
\overline{\tau}d\tau - \tau d\overline{\tau}) \Big\},$$ one checks that $d\theta^{\nu,\mu}_\tau = d\theta^{B^{\nu,\mu}_\tau},$ where $\theta^{B^{\nu,\mu}_\tau}:=iB^{\nu,\mu}_\tau (\bar z dz-zd\bar z)$. Therefore, there exists a function
$\widetilde{\varphi^{\nu,\mu}_\tau}: \C \to \C$ such that $\Im\widetilde{\varphi^{\nu,\mu}_\tau}= {Constant}$ and satisfying the first order partial differential equation
\begin{eqnarray} \frac{\partial \widetilde{\varphi^{\nu,\mu}_\tau}}{\partial \bar z} &=& -i\Big([\nu - B^{\nu,\mu}_\tau]z +\mu( \tau \frac{\partial \bar\tau}{\partial \bar z}- \bar \tau\frac{\partial \tau}{\partial\bar z}) \Big) \nonumber \label{Existence1stDiffEq}
\\ &=&  -i\mu\Big(( \tau \frac{\partial \bar\tau}{\partial \bar z}- \bar \tau\frac{\partial \tau}{\partial\bar z}) - (|\frac{\partial \tau}{\partial z}|^2- |\frac{\partial \tau}{\partial\bar z}|^2) z\Big) \label{Existence1stDiffEqP}
 .\end{eqnarray}
This completes the proof.
\end{proof}

\begin{remark} \label{RemReduction} The partial differential equation \eqref{Existence1stDiffEqP} %or \eqref{Existence1stDiffEq}
 satisfied by $\widetilde{\varphi^{\nu,\mu}_\tau}$
 can be reduced further to the following
 \begin{equation}
 \frac{\partial \psi^{\nu,\mu}_\tau}{\partial \bar z} =   \bar \tau \frac{\partial \tau}{\partial \bar z}
 \quad \mbox{ and } \quad \frac{\partial \psi^{\nu,\mu}_\tau}{\partial  z} =   \frac 1\mu (\nu-B^{\nu,\mu}_\tau)\bar z +\bar \tau \frac{\partial \tau}{\partial z} \label{1stDiffEq-Red}
   \end{equation}
  with
$\widetilde{\varphi^{\nu,\mu}_\tau}(z)= i\bigg([B^{\nu,\mu}_\tau -\nu]|z|^2 -\mu|\tau(z)|^2\bigg) + 2i\mu  \psi^{\nu,\mu}_\tau(z).$
\end{remark}

\section{Main result} %Proof of Theorem \ref{Thm2}}
 Let $\varphi^{\nu,\mu}_\tau$ be the real part of $\widetilde{\varphi^{\nu,\mu}_\tau}-\widetilde{\varphi^{\nu,\mu}_\tau}(0)$, where $\widetilde{\varphi^{\nu,\mu}_\tau}$ is a complex valued function on $\C$ as in Proposition \ref{1stdiffEq}. Define $\mathcal{W}^{\nu,\mu}_\tau $ to be the transformation given by
\begin{eqnarray}\label{GaugeTransformation}
[\mathcal{W}^{\nu,\mu}_\tau  (f)](z):= e^{i \varphi^{\nu,\mu}_\tau(z)}f(z).
\end{eqnarray}
We have
\begin{theorem} \label{Thm2}
  The image of $\mathcal{M}^{\nu,\mu}_{\tau}(\C)$ by the  transform  \eqref{GaugeTransformation}
  is the space of Landau $(\Gamma,\chi_\tau)$-automorphic functions. More exactly, we have
  $$\mathcal{W}^{\nu,\mu}_\tau (\mathcal{M}^{\nu,\mu}_{\tau}(\C) )=\set{F; ~ \mathcal{C}^\infty, \quad F(z+\gamma)= \chi_\tau(\gamma)j^{B^{\nu,\mu}_{\tau}}(\gamma, z)  F(z)}
  %:=\mathcal{F}^{B^{\nu,\mu}_{\tau}}_{\Gamma,\chi_{\tau}}
  ,$$
where $B^{\nu,\mu}_\tau = \nu +\mu(|\frac{\partial \tau}{\partial z}|^2- |\frac{\partial \tau}{\partial\bar z}|^2) \in \R$ and  $\chi_\tau$ is the pseudo-character defined on $\Gamma$ by
 $$\chi_\tau(\gamma)=\exp\big({2i\varphi^{\nu,\mu}_\tau(\gamma)-2i\mu  \Im \scal{\tau(0),\rho(\gamma)^{-1}\cdot 0}}\big).$$
 \end{theorem}

For the proof,  we begin with the following
 \begin{lemma}\label{LemTransMixed}
 The function $\widehat{\chi_\tau}$ defined on $\C \times \Gamma$ by
 $$\widehat{\chi_\tau}(z;\gamma):=  e^{i(\varphi^{\nu,\mu}_\tau(z+\gamma)- \varphi^{\nu,\mu}_\tau(z))} e^{-2i([B^{\nu,\mu}_{\tau}-\nu] \Im\scal{z,\gamma} +\mu  \Im \scal{\tau(z),\rho(\gamma)^{-1}\cdot 0})}$$ is independent of the variable $z$.
 \end{lemma}

 \begin{proof} Differentiation of $\widehat{\chi_\tau}(z;\gamma)$ w.r.t. variable $z$  gives
\begin{align} \frac{\partial \widehat{\chi_\tau}}{\partial z}   & =  i \bigg(
 \frac{\partial \varphi^{\nu,\mu}_\tau}{\partial z}(z+\gamma) - \frac{\partial\varphi^{\nu,\mu}_\tau}{\partial z}(z)\bigg) \widehat{\chi_\tau} \label{Indep1} \\ & \qquad
 - \bigg([B^{\nu,\mu}_{\tau}-\nu] \bar\gamma  + \mu  \bigg[ \overline{(\rho(\gamma^{-1})\cdot 0)} \frac{\partial\tau}{\partial z}(z) - (\rho(\gamma^{-1})\cdot 0) \frac{\partial\bar \tau}{\partial z}(z)\bigg]
 \bigg)\widehat{\chi_\tau}. \nonumber\end{align}
 In the other hand, using the equivariant condition $\tau(z+\gamma)= \rho(\gamma)\cdot \tau(z)$ and equation \eqref{Existence1stDiffEq}, one gets
 \begin{align} i\bigg(\frac{\partial \varphi^{\nu,\mu}_\tau}{\partial z}(z+\gamma)- \frac{\partial\varphi^{\nu,\mu}_\tau}{\partial z}(z)\bigg)
% = - (\overline{S^{\nu,\mu}_{\tau}}(z+\gamma) - \overline{(z+\gamma)})  - -  (\overline{S^{\nu,\mu}_{\tau}}(z) - \bar z)
& =   B^{\nu,\mu}_{\tau}\bar\gamma +  \overline{S^{\nu,\mu}_{\tau}}(z) - \overline{S^{\nu,\mu}_{\tau}}(z+\gamma) \nonumber\\
%& =   B^{\nu,\mu}_{\tau}\bar\gamma  -\nu \bar\gamma + \mu\bigg( \overline{\rho(\gamma)^{-1}\cdot 0} \frac{\partial\tau}{\partial z}(z) - \rho(\gamma)^{-1}\cdot 0 \frac{\partial\bar \tau}{\partial z}(z) \bigg)   \\
& =  [B^{\nu,\mu}_{\tau}-\nu]\bar\gamma   + \mu\bigg[ \overline{a_\gamma} \frac{\partial\tau}{\partial z}(z) - a_\gamma \frac{\partial\bar \tau}{\partial z}(z)\bigg] ,   \label{Indep2}
   \end{align}
 where we have set $a_\gamma=\rho(\gamma^{-1})\cdot 0$. Thus from \eqref{Indep1} and \eqref{Indep2}, we conclude that $\frac{\partial \widehat{\chi_\tau}}{\partial z}=0.$ Similarly, one gets also $\frac{\partial \widehat{\chi_\tau}}{\partial \bar z}=0$. This ends the proof of Lemma \ref{LemTransMixed}.
 \end{proof}

   \begin{proof}[Proof of Theorem \ref{Thm2}]
  We have to prove that $\mathcal{W}^{\nu,\mu}_\tau F$ belongs to
  $$\mathcal{F}^{B^{\nu,\mu}_{\tau}}_{\Gamma,\chi_{\tau}} : =\set{F ;\, \mathcal{C}^\infty, ~ F(z+\gamma)=
\chi_\tau(\gamma)j^{B^{\nu,\mu}_{\tau}}(\gamma, z)  F(z)}$$
   whenever  $F \in \mathcal{M}^{\nu,\mu}_{\tau}(\C)$, where $$\chi_\tau(\gamma):=\exp\big(2i\varphi^{\nu,\mu}_\tau(\gamma)-2i\mu  \Im \scal{\tau(0),\rho(\gamma)^{-1}\cdot 0}\big).$$
Indeed, we have
\begin{align*}
[\mathcal{W}^{\nu,\mu}_\tau F](z+\gamma)
& := e^{i\varphi^{\nu,\mu}_\tau(z+\gamma)} F(z+\gamma) \\
&  = e^{i\varphi^{\nu,\mu}_\tau(z+\gamma)} j^{\nu}(\gamma, z)j^{\mu}(\rho(\gamma), \tau(z)) F(z) \\
%& = e^{i(\varphi^{\nu,\mu}_\tau(z+\gamma)-\varphi^{\nu,\mu}_\tau(z))} j^{\nu}(\gamma, z)j^{\mu}(\rho(\gamma), \tau(z)) e^{i\varphi^{\nu,\mu}_\tau(z)} F(z) \\
&  = e^{i(\varphi^{\nu,\mu}_\tau(z+\gamma)-\varphi^{\nu,\mu}_\tau(z))} j^{\nu}(\gamma, z)j^{\mu}(\rho(\gamma), \tau(z)) [\mathcal{W}^{\nu,\mu}_\tau F](z) \\
&  = \widehat{\chi_\tau}(z;\gamma)  j^{-B^{\nu,\mu}_\tau}(\gamma, z) [\mathcal{W}^{\nu,\mu}_\tau F](z) .
\end{align*}
 Whence by Lemma  \ref{LemTransMixed}, we see that  $\widehat{\chi_\tau}(z;\gamma) = \widehat{\chi_\tau}(0;\gamma) =: \chi_\tau(\gamma) $ and therefore
 $$ [\mathcal{W}^{\nu,\mu}_\tau F](z+\gamma) = \chi_\tau(\gamma)  j^{-B^{\nu,\mu}_\tau}(\gamma, z) [\mathcal{W}^{\nu,\mu}_\tau F](z).$$
  The proof is completed
 \end{proof}

 \begin{corollary} \label{CorRDQ}
 The function $\chi_\tau(\gamma)=\exp\big(2i\varphi^{\nu,\mu}_\tau(\gamma)-2i\mu  \Im \scal{\tau(0),\rho(\gamma)^{-1}\cdot 0}\big)$
 satisfies the pseudo-character property
 $$\chi_\tau(\gamma+\gamma')=e^{2i B^{\nu,\mu}_\tau \Im \scal{\gamma,\gamma'} } \chi_\tau(\gamma)\chi_\tau(\gamma')$$
 if and only if  $(1/\pi)\phi^{\nu,\mu}_\rho$ in \eqref{PhaseFactor} takes integer values on $\Gamma\times \Gamma$.
 \end{corollary}

%\section{Concluding remarks}

 {\bf\it Acknowledgements.} The author is indebted to Professor A. Intissar for valuable discussions and encouragement.

\end{document}